# Evaluating the potential of HIV self-testing to reduce HIV incidence in EHE districts: a modeling study

*Alex Viguerie, PhD*[a,*]

*Dipartimento di Scienze Pure ed Applicate*
*Università degli Studi di Urbino Carlo Bo*
*Via Aurelio Saffi 2*
*Urbino, PU 61029 (Italy)*

alexander.viguerie@uniurb.it

**Evaluating the potential of HIV self-testing to reduce HIV incidence in EHE districts: a modeling study**

**Abstract**

**Background**
High HIV transmission persists in many U.S. jurisdictions despite prevention efforts. HIV self-testing offers a means to overcome barriers associated with routine laboratory-based testing but carries a risk of increasing incidence if replacement effects reduce overall test sensitivity.

**Methods**
A linearized four-compartment HIV transmission model was applied to 38 Ending the HIV Epidemic (EHE) priority jurisdictions. Both the percentage of self-tests and the overall testing rate were varied to quantify 10-year changes in HIV incidence. A threshold testing level was defined to counterbalance potential negative effects from reduced self-test sensitivity.

**Results**
Substantial heterogeneity emerged across districts. Incidence reductions exceeded 5% in some areas, while others saw only minor effects. Jurisdictions with higher baseline testing displayed an elevated risk of increased incidence from substitution of laboratory-based testing with self-tests. In contrast, a derived *Awareness Reproduction Number*—capturing transmissions attributable to undiagnosed infection—strongly correlated with the magnitude of possible incidence declines.

**Conclusions**
Local epidemiological context is pivotal in determining the risks and benefits of HIV self-testing. Jurisdictions with robust testing systems may face a greater likelihood of inadvertently raising incidence, whereas those with a high number of individuals stand to achieve notable transmission reductions. Tailoring self-testing strategies to each district's baseline conditions can maximize public health benefits while minimizing unintended consequences.

## Introduction

Despite notable progress in HIV prevention and treatment, many communities in the United States continue to experience substantial transmission rates and suboptimal testing uptake. One promising approach to expanding testing access is HIV self-testing, which can help overcome barriers associated with stigma, inconvenient testing hours or locations, and the logistical challenges of routine laboratory-based testing. In principle, HIV self-testing can bridge gaps in the HIV care continuum by reaching those individuals who have never or rarely tested, facilitating earlier diagnosis, and potentially reducing incidence—provided that positive self-test results are promptly confirmed and linked to care [1], [2], [3], [4], [5], [6], [7], [8].

Recent modeling work suggests that increased HIV self-testing in the US MSM population may reduce HIV incidence by as much as 10% over a ten-year horizon [8]. Furthermore, although HIV self-tests are less sensitive and have a longer window period compared to laboratory-based testing [9], [10], [11], [12], the risk of increased incidence and/or reduced awareness of status due to replacement of gold-standard laboratory tests with HIV self-tests was found to be minimal [8]. These findings differed somewhat from previous modeling studies, which found similar potential HIV incidence reductions; but with a notably greater risk of potential negative effects, including increased HIV incidence owing to widespread HIV self-testing [2], [4].

These conflicting findings suggest a need for further modeling work on this question. Notably, the analysis [8] considered the entire US MSM population, while [2], [4] considered specific regional populations. Given the widespread differences in HIV prevalence, incidence, awareness, and viral suppression across US jurisdictions, the differing conclusions in modeling studies are not surprising. Indeed, such findings suggest that the effect of increased HIV self-testing may vary significantly by jurisdiction, in terms of both potential incidence reductions, as well as the possible risk of increased transmission.

The present study attempts to further explore the heterogeneous impacts of expanded HIV self-testing across jurisdictions by applying the modeling framework introduced in [8] across 38 Ending the HIV Epidemic (EHE) [13], [14], [15] priority districts with differing HIV profiles. For each district, we quantify both the potential benefits of expanded HIV self-testing, in terms of reduced HIV incidence and increased awareness, and the potential risks, including the potential for increased HIV incidence.

Our findings show a decoupling between the potential risks and benefits of expanded HIV testing. We find that the risk of increased HIV incidence in a jurisdiction is primarily driven by the baseline jurisdictional awareness of HIV status. In contrast, we find that the potential reductions in HIV incidence are driven by a jurisdiction-dependent quantity, termed the *Awareness Reproduction Number*. We provide both a mathematical derivation of this quantity, and establish its relationship between HIV self-testing and HIV incidence both empirically. Using this risk/benefit framework, we establish criteria for identifying jurisdictions offering a combination of low risk and high benefit, in which expanded HIV self-testing interventions may result in significant HIV reductions at minimal risk.

## Methods

*Compartmental Model*

We consider the four-compartment model introduced in [8], and refer the reader to that work for the detailed analysis and explanation of the model.

The four compartments correspond to the disease stages: acute HIV infection $a$, unaware chronic HIV infection $u$, unaware AIDS infection $s$ and diagnosed HIV infection $d$. The full system of equations reads as follows:

$$\begin{aligned}
\dot{a} &= \lambda_a a + \lambda_u u + \lambda_s s + \lambda_d d - \left(\sigma_{a\to u} + \widetilde{\phi_a} + \mu_a\right)a \\
\dot{u} &= \sigma_{a\to u}\, a - \left(\sigma_{u\to s} + \widetilde{\phi_u} + \mu_u\right)u \qquad (1) \\
\dot{s} &= \sigma_{u\to s} u - \left(\widetilde{\phi_s} + \mu_s\right)s \\
\dot{d} &= \widetilde{\phi_a} a + \widetilde{\phi_u} u + \widetilde{\phi_s} s - \mu_d d,
\end{aligned}$$

together with appropriate initial conditions. Table 1 reports the values, units and sources for the quantities these considerations in mind, the description of each parameter and the associated units are provided in Table 1. A full explanation of how the model is parameterized from surveillance data is provided in Supplement B.

The model *detection rates* are given by:

$$\widetilde{\phi_k} = \kappa_k^{self} \gamma_k \left[\frac{1}{\left(\frac{1}{(1+\chi_k)\phi_k} + t_k^{self\to diag.}\right)}\right] + \kappa_k^{care}(1-\gamma_k)(1+\chi_k),$$

where $\phi_k$ is the baseline jurisdictional HIV testing rate and $\kappa_k^{self}$ and $\kappa_k^{care}$ are sensitivities for HIV self-tests and non self-tests, respectively.

HIV self-testing interventions are modeled through by varying the parameters $\gamma_k$ and $\chi_k$, which refers to the percentage of HIV tests that are self-tests, and the percent increase in overall HIV testing from the jurisdictional baseline level. At baseline, both are set to zero. Finally, the parameter $t_k^{self\to diag.}$ refers to the time to formal HIV diagnosis after a positive HIV self-test.

In Supplement A, we demonstrate that the above model can be regarded as the *linearized infection subsystem* of a nonlinear model of HIV transmission. Additionally, in Supplement A, we rigorously establish that the approximation error arising from the model linearization remains small for realistic parameter choices and medium-range (10-20 year) time intervals.

*Simulation study design*

To evaluate the possible effects of expanded HIV self-testing, apply the model to a range of 38 jurisdictions in the United States, each considered an EHE priority district. For each district, we perform a series of simulations, varying both the percentage of HIV self-tests among total tests from 0% to 100%, and the percent increase in overall HIV testing from 0% to 100%. Each simulation considers a ten-year intervention period, and the primary outcome of interest is the percent change in cumulative 10-year HIV incidence compared to the baseline scenario (no self-testing and baseline testing rates throughout the 10-year simulation period). A comparison of the baseline model configuration for each jurisdiction, with surveillance data is provided in Supplement B.

To quantify the risk associated with increased HIV self-testing, we calculate *threshold testing level* $\chi_\gamma$ for a considered proportion of HIV self-tests $\gamma$. For a given proportion of HIV self-tests $\gamma$, the threshold testing level $\chi_\gamma$ is defined as the increase in HIV testing rate necessary to ensure that HIV incidence is decreased, offsetting any possible negative effects from reduced HIV test sensitivity and/or increased window period. Note that higher threshold levels indicates a higher-risk jurisdiction.

To quantify the potential benefits of HIV self testing, we calculate, for each jurisdiction, the mean percent reduction from baseline (ie, no changes in testing rates and no self-testing) across all considered HIV self-testing scenarios in the jurisdiction.

To connect our modeling analysis to actionable policy, we further analyze the relationship between both risk (as quantified by $\chi_\gamma$) and benefits (as quantified by mean incidence reduction) to baseline-level indicators for each jurisdiction. These indicators include:

a. Baseline viral reproduction number $R_t$ (definition in Supplement A). Intuitively, this can be regarded as the expected number of transmissions per PWH with newly acquired infection.
b. Baseline awareness reproduction number $R_{Awr}$ (definition in Supplement A). Intuitively, this can be regarded as the expected number of transmissions per PWH with newly acquired infection, attributable to unawareness of HIV status.
c. Baseline transmission rate per month $\overline{\lambda}$ (definition in Supplement B)
d. Baseline testing rate per month $\overline{\phi}$ (definition in Supplement B

**Results**

In Table 2, we provide the results for each jurisdiction. Our outputs of interest are:

- Percent reduction in 10-year HIV incidence across all self-testing scenarios. This quantifies the possible benefits of increased self-testing.
- Threshold testing levels for $\chi_\gamma$ for $\gamma = 0.25,\ 0.5, 0.75, 1.0$. These are the necessary increases in overall testing rate necessary to offset reduced sensitivity when self-tests comprise, respectively, 25%, 50%, 75%, 100%, of all HIV testing. These indicators quantify the inherent risk of replacement effects due to increased self-testing.

We observe wide variation in potential incidence reductions, with King County, WA showing a 7.3% reduction in HIV incidence, and New York County, NY showing a 1.5% reduction in HIV incidence. Over all jurisdictions, the mean percent reduction in incidence was 4.0%, and the median was 3.9%.

In Figure 2, we plot the association between different baseline indicators in each jurisdiction and the percent incidence reduction. We see very little association between incidence reduction due to HIV self-testing and the baseline incidence and testing rates, or the baseline effective reproduction number $R_t$. However, the baseline Awareness Reproduction Number $R_{Awr}$ (full derivation and definition in Supplement A), shows a near-perfect association with potential incidence reduction.

The threshold testing levels also varied across jurisdictions, with variation becoming more pronounced as the proportion of self-testing increased. San Francisco County, CA consistently showed the largest necessary increases in HIV testing, while San Diego County, CA showed the lowest necessary increases. At $\gamma = 0.25$, we find very little difference across jurisdictions, with the threshold testing level between 3.8% and 5.4% for all jurisdictions. This increases as the proportion of self-tests increases, and we see a nearly 10 percentage- point spread at $\gamma = 1.0$ between the highest (27.3%) and lowest (17.4%) threshold testing levels.

In Figure 3, we plot the association between baseline jurisdictional indicators and threshold testing levels. Little association is observed between threshold testing levels and baseline viral reproduction number $R_t$ and baseline HIV incidence rate. However, a strong association is observed between threshold testing levels and baseline HIV testing rate.

**Table 1:** Model parameters and values.

| Parameter | Name | Value | Source |
|---|---|---|---|
| $\sigma_{a \to u}$ | Rate of movement from acute to chronic infection | 1/60 days | [17] |
| $\sigma_{u \to s}$ | Rate of movement from chronic infection to AIDS | 1/11.8 years | [18] |
| $t_s^{self \to diag.}$ | Delay in diagnosis after positive self-test, AIDS | 1/30 days | Assumed |
| $t_{a,u}^{self \to diag.}$ | Delay in diagnosis after positive self-test, unaware acute, unaware chronic | 1/90 days | Assumed |
| $\phi_{a,u}$ | Baseline testing rate for PWH with acute and chronic infection | .0118 / months | Calculated using Surveillance data, 2017-19 [16]; see Supplement B |
| $\phi_s$ | Testing rate for PWH with AIDS (late-stage) | .0481/months | Calculated using Surveillance data, 2017-19 [16]; set as $\nu_s \phi_u$see Supplement B |
| $\nu_s$ | Multiplier for increase in testing from AIDS compared to acute/chronic infection | 4.08 | [19] |
| $\kappa_a^{care}$ | Sensitivity of laboratory test to acute infection | .83 | [17] |
| $\kappa_a^{self}$ | Sensitivity of self-test to acute infection | 0.0 | [10] |
| $\kappa_{a,u}^{care}$ | Sensitivity of laboratory test to chronic infection | 1.0 | [20] |
| $\kappa_{u,s}^{self}$ | Sensitivity of self-test to chronic infection | .92 | [10] |
| $\mu_a$ | Mortality rate, PWH with acute infection | .0069/years | Surveillance data, 2017-19 [16]; see Supplement B |
| $\mu_u, \mu_{noCare}, \mu_{ART}$ | Mortality rate, PWH with chronic infection, diagnosed not in care, on ART but not VLS | .0174/years | Surveillance data, 2017-19 [16]; set as$\beta_u \mu_a$, $\beta_{noCare} \mu_a$, $\beta_{ART} \mu_a$; see Supplement B |
| $\beta_u, \beta_{noCare}, \beta_{ART}$ | Factor change in mortality for PWH with chronic infection, diagnosed but not in care, and ART but not VLS compared to acute | 2.538 | [21], assumed equal. |
| $\mu_s$ | Mortality rate, PWH with AIDS, not diagnosed | .046/years | [21], see Supplement B, set as $\beta_s \mu_a$ |
| $\beta_s$ | Factor change in mortality for PWH with undiagnosed AIDS compared to PWH with acute infection | 6.172 | [21] |
| $\mu_d$ | Mortality rate, PWH with diagnosed infection | .0086/years | Surveillance data, 2017-19 [16]; see Supplement B |
| $\mu_{VLS}$ | Mortality rate, PWH who are VLS | .0043/years | Surveillance data, 2017-19 [16]; set as $\beta_{VLS} \mu_a$, see Supplement B |
| $\beta_{VLS}$ | Factor change in mortality for PWH who are VLS compared to acute infection | 0.6346 | [21] |

## Discussion

Our results show that the impact of expanded HIV self-testing, in terms of its potential risks and benefits, varies by jurisdiction, in some cases significantly. The variation is most pronounced in terms of incidence reduction, with several jurisdictions showing quite minimal benefit (less than 2% reduction in HIV incidence over ten years). In contrast, several other jurisdictions show greater potential for HIV incidence reduction, as much as 6% or 7%. While the risk profiles vary less than possible benefits, at higher levels of self-test proportions, we observe noticeable differences, with necessary testing increases differing by as much as 10% across jurisdictions.

Importantly, our results suggest that the jurisdictional risks and benefits appear independent. Indeed, among low-benefit (<2% incidence reduction) districts, there are districts with high (e..g. San Francisco County, CA) and low (e.g. Baltimore County, MD) threshold-testing levels. Similarly, our analysis showed high-benefit jurisdictions with lower (e.g. Travis County, TX) and higher (e.g. Gwinett County, GA) threshold testing levels.

To better understand the drivers of each of the risk and benefit components, we examined the relationship between jurisdictional risk (as measured by threshold testing levels) and benefits (as measured by percent 10-year incidence reduction from baseline) and baseline jurisdictional quantities, which are either known or possible to estimate from surveillance data. Jurisdictional risk is primarily determined by the baseline jurisdiction testing rate, with other indicators showing little or no association with threshold testing levels.

The primary driver of the potential benefits is less straightforward, as there is no clear relationship between potential incidence reduction, baseline testing rate, baseline incidence rate, or effective viral reproduction number. However, through a more detailed mathematical analysis (shown in Supplement A), we found that the derived *awareness reproduction number* bears a near-perfect relationship between the observed variation in potential incidence reduction across jurisdictions.

This analysis has several key limitations. As discussed in [8], while the modeling approach here has many attractive aspects (including being directly parameterizable from surveillance data), it may not be well-suited for districts in which the population transmission profile is highly irregular or nonlinear. Such dynamics are commonly observed in e.g., outbreaks among injection drug users [22]. Additionally, the model assumes that self-testing rates lead to increased self-testing among undiagnosed PWH; while this is of course the goal of HIV self-testing, it may not always be the case in practice [5]. Furthermore, while evidence has shown the potential of HIV self-testing to increase testing among key groups, particularly in programmatic settings [3], [6], [23], the current modeling study does not provide any guidelines or insights into the design and/or implementation of such programs.

Nonetheless, the analysis shows several important findings with implications for broader HIV self-testing policy. As suggested by previous modeling studies, which gave mixed results [2], [4], [8], both the possible benefits from expanded HIV self-testing, and the risk of HIV self-testing increasing HIV incidence, depend on jurisdiction-specific factors, and are largely independent of each other. In practice, districts and jurisdictions with high baseline testing rates show a higher risk of increased HIV transmission due to increased self-testing. The primary driver of potential incidence reduction is not as straightforward, however, our current analysis shows that an indicator can be derived, based on HIV surveillance data, that aggregates several pieces of relevant information and quantifies HIV transmission due to unawareness of HIV status. Therefore, this study provides policymakers with a clear set of criteria for identifying priority districts in which HIV self-testing may deliver large benefits at low costs.

**Table 2:** Simulation results for each jurisdiction. $\overline{\lambda}$ is the baseline transmission rate, $\overline{\phi}$ the baseline testing rate, $R_t$ the viral reproduction number, $R_{Awr}$ the awareness reproduction number, % inc. red. refers to the percentage reduction in 10-year HIV incidence due to increased self-testing. The $\chi_\gamma$ refer to threshold testing levels – necessary increase in testing rate to offset reduced sensitivity for a given proportion of self-tests $\gamma$.

| **Jurisdiction** | $\overline{\lambda}$ | $\overline{\phi}$ | $R_t$ | $R_{Awr}$ | **% inc.red.** | $\chi_{0.25}$ | $\chi_{0.5}$ | $\chi_{0.75}$ | $\chi_{1.0}$ |
|---|---|---|---|---|---|---|---|---|---|
| Alameda County, CA | 0.028 | 0.016 | 2.539 | 0.253 | *5.0%* | 4.0% | 8.5% | 13.5% | 18.9% |
| Los Angeles County, CA | 0.027 | 0.018 | 2.490 | 0.169 | *3.3%* | 4.2% | 8.9% | 14.2% | 19.8% |
| Riverside County, CA | 0.028 | 0.016 | 1.609 | 0.294 | *5.6%* | 4.1% | 8.7% | 13.6% | 18.9% |
| Sacramento County, CA | 0.048 | 0.014 | 2.887 | 0.338 | *6.7%* | 4.0% | 8.3% | 13.2% | 18.2% |
| San Bernardino County, CA | 0.051 | 0.015 | 3.349 | 0.230 | *4.7%* | 4.1% | 8.6% | 13.8% | 19.1% |
| San Diego County, CA | 0.029 | 0.014 | 2.819 | 0.175 | *3.4%* | 3.8% | 8.0% | 12.5% | 17.4% |
| San Francisco County, CA | 0.012 | 0.028 | 1.674 | 0.128 | *2.2%* | 5.4% | 11.6% | 19.0% | 27.3% |
| Broward County, FL | 0.027 | 0.017 | 1.825 | 0.199 | *3.9%* | 4.2% | 8.8% | 14.0% | 19.5% |
| Duval County, FL | 0.033 | 0.016 | 1.767 | 0.173 | *3.4%* | 4.1% | 8.7% | 13.7% | 19.1% |
| Hillsborough County, FL | 0.033 | 0.015 | 2.005 | 0.238 | *4.7%* | 4.0% | 8.4% | 13.4% | 18.5% |
| Miami-Dade County, FL | 0.038 | 0.018 | 2.562 | 0.193 | *3.9%* | 4.4% | 9.3% | 14.8% | 20.8% |
| Palm Beach County, FL | 0.028 | 0.015 | 1.523 | 0.161 | *3.1%* | 4.0% | 8.4% | 13.3% | 18.4% |
| Pinellas County, FL | 0.027 | 0.017 | 1.332 | 0.220 | *4.2%* | 4.2% | 8.8% | 13.9% | 19.4% |
| Cobb County, GA | 0.050 | 0.015 | 5.256 | 0.265 | *5.6%* | 4.1% | 8.5% | 13.7% | 18.9% |
| Dekalb County, GA | 0.037 | 0.016 | 3.463 | 0.213 | *4.3%* | 4.0% | 8.5% | 13.5% | 18.8% |
| Fulton County, GA | 0.034 | 0.016 | 3.113 | 0.187 | *3.8%* | 4.0% | 8.5% | 13.5% | 18.8% |
| Gwinnett County, GA | 0.048 | 0.017 | 5.727 | 0.267 | *5.7%* | 4.2% | 9.0% | 14.3% | 19.9% |
| Marion County, IN | 0.045 | 0.016 | 3.575 | 0.258 | *5.4%* | 4.2% | 8.9% | 14.1% | 19.5% |
| East Baton Rouge Parish, LA | 0.039 | 0.016 | 1.798 | 0.272 | *5.4%* | 4.2% | 8.7% | 13.9% | 19.4% |
| Orleans Parish, LA | 0.027 | 0.019 | 1.747 | 0.183 | *3.6%* | 4.4% | 9.4% | 14.9% | 21.0% |
| Montgomery County, MD | 0.022 | 0.020 | 3.894 | 0.104 | *2.0%* | 4.6% | 9.7% | 15.6% | 22.0% |
| Prince George's County, MD | 0.025 | 0.017 | 3.898 | 0.146 | *2.9%* | 4.1% | 8.8% | 13.9% | 19.4% |
| Baltimore City, MD | 0.016 | 0.017 | 0.889 | 0.093 | *1.7%* | 4.1% | 8.8% | 14.1% | 19.5% |
| Wayne County, MI | 0.036 | 0.018 | 2.035 | 0.243 | *4.9%* | 4.3% | 9.0% | 14.5% | 20.4% |
| Bronx County, NY | 0.014 | 0.019 | 0.884 | 0.097 | *1.7%* | 4.3% | 9.2% | 14.7% | 20.6% |
| Kings County, NY | 0.016 | 0.018 | 1.118 | 0.109 | *2.0%* | 4.3% | 9.0% | 14.3% | 20.2% |
| New York County, NY | 0.011 | 0.016 | 1.066 | 0.079 | *1.5%* | 4.0% | 8.4% | 13.4% | 18.6% |
| Queens County, NY | 0.017 | 0.018 | 1.923 | 0.117 | *2.3%* | 4.2% | 8.9% | 14.2% | 19.7% |
| Cuyahoga County, OH | 0.024 | 0.015 | 1.801 | 0.130 | *2.5%* | 4.0% | 8.3% | 13.2% | 18.2% |
| Franklin County, OH | 0.036 | 0.017 | 2.811 | 0.219 | *4.4%* | 4.2% | 8.9% | 14.2% | 19.6% |
| Hamilton County, OH | 0.038 | 0.018 | 2.544 | 0.133 | *2.6%* | 4.4% | 9.4% | 15.0% | 21.0% |
| Philadelphia County, PA | 0.023 | 0.019 | 1.548 | 0.137 | *2.6%* | 4.4% | 9.5% | 15.3% | 21.5% |
| Shelby County, TN | 0.036 | 0.017 | 2.113 | 0.208 | *4.2%* | 4.3% | 9.0% | 14.4% | 20.0% |
| Dallas County, TX | 0.037 | 0.016 | 3.957 | 0.200 | *4.1%* | 4.0% | 8.5% | 13.4% | 18.7% |
| Harris County, TX | 0.038 | 0.016 | 2.664 | 0.192 | *3.9%* | 4.1% | 8.7% | 13.8% | 19.2% |
| Tarrant County, TX | 0.049 | 0.016 | 3.149 | 0.276 | *5.7%* | 4.2% | 8.7% | 13.9% | 19.4% |
| Travis County, TX | 0.037 | 0.014 | 2.995 | 0.323 | *6.3%* | 3.9% | 8.1% | 12.7% | 17.6% |
| King County, WA | 0.029 | 0.017 | 2.212 | 0.370 | *7.3%* | 4.2% | 8.6% | 13.9% | 19.3% |

**Figure 2:** A series of plots comparing reduction in incidence in response to increased HIV self-testing and different baseline jurisdictional indicators. We note that the derived Awareness Reproduction Number $R_{Awr}$ (see Supplement A for derivation) provides the strongest quantifier of potential benefit from HIV self-testing. Overall incidence rate, testing rate, and viral reproduction number $R_t$ appear uncorrelated.

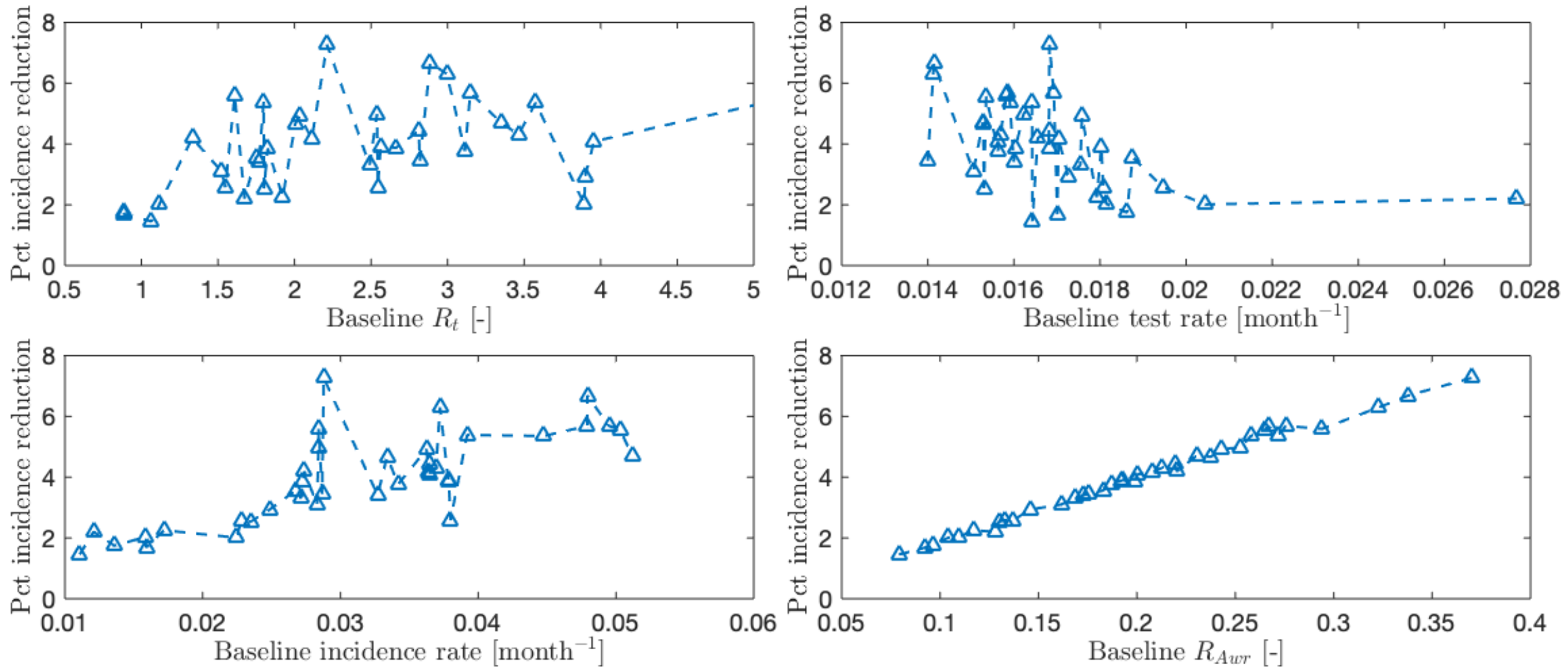

**Figure 3:** A series of plots comparing the threshold testing levels, which quantify the increase in testing rate necessary to ensure HIV self-testing results in reduced HIV incidence, to jurisdictional baseline indicators. The only strong association is with the baseline jurisdictional testing rate. Overall incidence rate, and viral reproduction number $R_t$ appear uncorrelated.

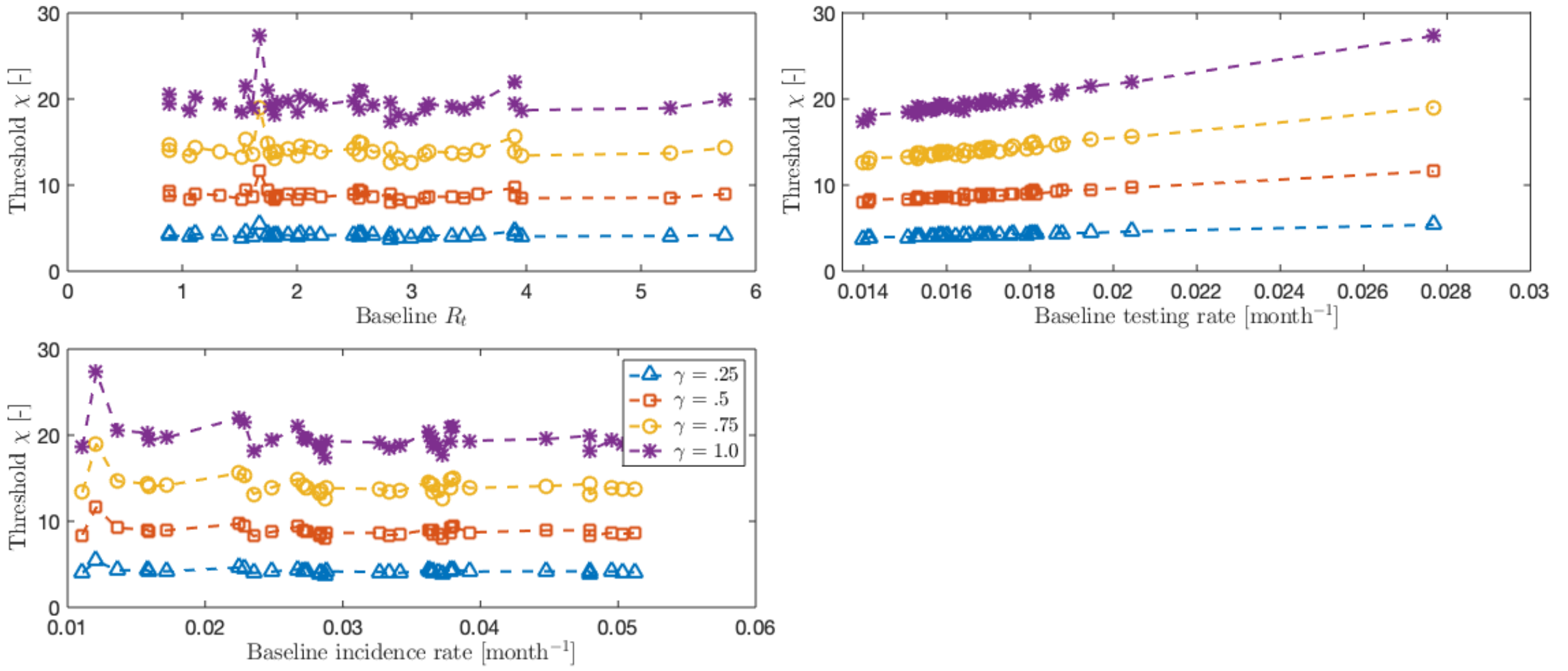

## Supplement A

## Mathematical analysis

Consider the following nonlinear model, related to (1) in the main text:

$$\begin{aligned}
\dot{e} &= \Lambda - \mu_e e - \left(\frac{\widetilde{\lambda_a} a}{n} + \frac{\widetilde{\lambda_u} u}{n} + \frac{\widetilde{\lambda_s} s}{n} + \frac{\widetilde{\lambda_d} d}{n}\right) e \\
\dot{a} &= \left(\frac{\widetilde{\lambda_a} a}{n} + \frac{\widetilde{\lambda_u} u}{n} + \frac{\widetilde{\lambda_s} s}{n} + \frac{\widetilde{\lambda_d} d}{n}\right) e - \left(\sigma_{a\to u} + \widetilde{\phi_a} + \mu_a\right) a \\
\dot{u} &= \sigma_{a\to u}\, a - \left(\sigma_{u\to s} + \widetilde{\phi_u} + \mu_u\right) u \qquad (A1) \\
\dot{s} &= \sigma_{u\to s} u - \left(\widetilde{\phi_s} + \mu_s\right) s \\
\dot{d} &= \widetilde{\phi_a} a + \widetilde{\phi_u} u + \widetilde{\phi_s} s - \mu_d d.
\end{aligned}$$

The parameters are identical to those in (1) in the main text, with the addition of the compartment $e$ (for *eligible*) denotes the eligible (or susceptible) population, $n = e + a + u + s + d$ is the total living population, $\Lambda$ (units Persons/Time) is the recruitment rate, and $\mu_e$ (units 1/Time) is the eligible mortality rate. $\widetilde{\lambda_i}$, $i = a, u, s, d$ have the same function and units as $\lambda_i$ in main text (1) (though not necessarily the same values).

We further denote the *eligible (susceptible) population fraction* as:

$$\Sigma = \frac{e}{n}$$

***Proposition 1:*** The system (1) in the main text is the *linearized infection subsystem* of (A1) linearized about $\Sigma_0 = \Sigma(t_0)$.

*Proof*: We follow the procedure outlined in [24], [25], [26]. Note that, in (A1), the infection can be transferred by all compartments besides $e$. Define:

$$\boldsymbol{x} = \begin{pmatrix} a \\ u \\ s \\ d \end{pmatrix}$$

and write (A1) as:

$$\dot{\boldsymbol{x}} = f(\boldsymbol{x}) - v(\boldsymbol{x}),$$

where:

$$f(\boldsymbol{x}) = \begin{pmatrix} \left(\widetilde{\lambda_a} a + \widetilde{\lambda_u} u + \widetilde{\lambda_s} s + \widetilde{\lambda_d} d\right)\Sigma \\ 0 \\ 0 \\ 0 \end{pmatrix}, v(\boldsymbol{x}) = \begin{pmatrix} \left(\sigma_{a\to u} + \widetilde{\phi_a} + \mu_a\right) a \\ -\sigma_{a\to u}\, a + \left(\sigma_{u\to s} + \widetilde{\phi_u} + \mu_u\right) u \\ -\sigma_{u\to s}\, u + \left(\widetilde{\phi_s} + \mu_s\right) s \\ \widetilde{-\phi_a} a - \widetilde{\phi_u} u - \widetilde{\phi_s} s + \mu_d d \end{pmatrix}.$$

Next, compute the *Jacobian matrices $F$ and $V$ of $f$ and $v$,* respectively, evaluated at $\Sigma = \Sigma_0$:

$$F_{i,j} = \frac{df(\boldsymbol{x_i})}{d\boldsymbol{x_j}}|_{\Sigma=\Sigma_0}, V_{i,j} = \frac{dv(\boldsymbol{x_i})}{d\boldsymbol{x_j}}|_{\Sigma=\Sigma_0}.$$

These yield:

$$F=\begin{pmatrix}\widetilde{\lambda_a}\times\Sigma_0 & \widetilde{\lambda_u}\times\Sigma_0 & \widetilde{\lambda_s}\times\Sigma_0 & \widetilde{\lambda_d}\times\Sigma_0\\ 0&0&0&0\\ 0&0&0&0\\ 0&0&0&0\end{pmatrix},\quad (A2)$$

$$V=\begin{pmatrix}\sigma_{a\to u}+\widetilde{\phi_a}+\mu_a & 0 & 0 & 0\\ -\sigma_{a\to u} & \sigma_{u\to s}+\widetilde{\phi_u}+\mu_u & 0 & 0\\ 0 & -\sigma_{u\to s} & \widetilde{\phi_s}+\mu_s & 0\\ -\widetilde{\phi_a} & -\widetilde{\phi_u} & -\widetilde{\phi_s} & \mu_d\end{pmatrix}.\,(A3)$$

The linearized infection subsystem about $\Sigma_0$ then reads:

$$\dot{\boldsymbol{x}}=(F-V)\boldsymbol{x}.$$

By letting $\lambda_i=\widetilde{\lambda_i}\times\Sigma_0$ for $i=a,u,s,d$, one obtains:

$$\begin{aligned}\dot{a}&=\lambda_a a+\lambda_u u+\lambda_s s+\lambda_d d-(\sigma_{a\to u}+\widetilde{\phi_a}+\mu_a)a\\ \dot{u}&=\sigma_{a\to u}a-(\sigma_{u\to s}+\widetilde{\phi_u}+\mu_u)u\\ \dot{s}&=\sigma_{u\to s}u-(\widetilde{\phi_s}+\mu_s)s\\ \dot{d}&=\widetilde{\phi_a}a+\widetilde{\phi_u}u+\widetilde{\phi_s}s-\mu_d d,\end{aligned}\qquad (A4)$$

which was to be shown.

***Proposition 2:*** The effective reproduction number $R_t$ for $\Sigma=\Sigma_0$ of (A1) and (A4) reads:

$$R_t=\frac{\lambda_a}{(\sigma_{a\to u}+\widetilde{\phi_a}+\mu_a)}+\frac{\lambda_u\sigma_{a\to u}}{(\sigma_{a\to u}+\widetilde{\phi_a}+\mu_a)(\sigma_{u\to s}+\widetilde{\phi_u}+\mu_u)}+\frac{\lambda_s\sigma_{a\to u}\sigma_{u\to s}}{(\sigma_{a\to u}+\widetilde{\phi_a}+\mu_a)(\sigma_{u\to s}+\widetilde{\phi_u}+\mu_u)(\widetilde{\phi_s}+\mu_s)}+\frac{\lambda_d}{\mu_d}\left[\frac{\widetilde{\phi_a}(\sigma_{u\to s}+\widetilde{\phi_u}+\mu_u)(\widetilde{\phi_s}+\mu_s)+\widetilde{\phi_u}\sigma_{a\to u}(\widetilde{\phi_s}+\mu_s)+\widetilde{\phi_s}\sigma_{a\to u}\sigma_{u\to s}}{(\sigma_{a\to u}+\widetilde{\phi_a}+\mu_a)(\sigma_{u\to s}+\widetilde{\phi_u}+\mu_u)(\widetilde{\phi_s}+\mu_s)}\right]$$

*Proof:*

$R_t$ is given by the largest eigenvalue of the *next generation matrix* $NV^{-1}$, where $N$ and $V$ are the matrices defining the linearized infection subsystem [24], [25], [26], [27]. As shown previously, for the system (A1), $N$ and $V$ are respectively given by (A2) and (A3).

Furthermore:

$$V^{-1}=\begin{pmatrix}\frac{1}{V_{1,1}} & 0 & 0 & 0\\ \frac{\sigma_{a\to u}}{V_{1,1}V_{2,2}} & \frac{1}{V_{2,2}} & 0 & 0\\ \frac{\sigma_{a\to u}\sigma_{u\to s}}{V_{1,1}V_{2,2}V_{3,3}} & \frac{\sigma_{u\to s}}{V_{2,2}V_{3,3}} & \frac{1}{V_{3,3}} & 0\\ \frac{\widetilde{\phi_a}}{V_{1,1}\mu_d}+\frac{\widetilde{\phi_u}\sigma_{a\to u}}{V_{1,1}V_{2,2}\mu_d}+\frac{\widetilde{\phi_s}\sigma_{a\to u}\sigma_{u\to s}}{V_{1,1}V_{2,2}V_{3,3}\mu_d} & \frac{\widetilde{\phi_u}}{V_{2,2}\mu_d}+\frac{\sigma_{u\to s}\widetilde{\phi_s}}{V_{2,2}V_{3,3}\mu_d} & \frac{\widetilde{\phi_s}}{V_{3,3}\mu_d} & \frac{1}{\mu_d}\end{pmatrix}$$

, where $V_{i,j}$ indicates the $i,j$-th entry of $V$.

Denoting the $i$-th column of $V^{-1}$ as $v_i$, and letting $n = (\lambda_a, \lambda_u, \lambda_s, \lambda_d)^T$, we observe that:

$$NV^{-1} = \begin{pmatrix} n^T v_1 & n^T v_2 & n^T v_3 & n^T v_4 \\ 0 & 0 & 0 & 0 \\ 0 & 0 & 0 & 0 \\ 0 & 0 & 0 & 0 \end{pmatrix},$$

$R_t$ is then obtained as the *spectral radius* of $NV^{-1}$, $\rho(NV^{-1})$, the magnitude of the largest eigenvalue of $NV^{-1}$ [24], [25], [26], [27].

Since that $NV^{-1}$ has only one non-zero row, $NV^{-1}$ has rank one (and therefore only one nonzero eigenvalue [28]). The only eigenvector corresponding to a nonzero eigenvalue is $(1,0,0,0)^T$ with eigenvalue $n^T v_1$. Thus:

$$|n^T v_1| = \rho(NV^{-1}) = R_t.$$

***Definition:*** We define the *Awareness Reproduction Number* as:

$$\begin{aligned} R_{Awr} &= R_t - \frac{\lambda_d}{\mu_d} \\ &= \frac{\lambda_a}{(\sigma_{a\to u} + \widetilde{\phi_a} + \mu_a)} + \frac{\lambda_u \sigma_{a\to u}}{(\sigma_{a\to u} + \widetilde{\phi_a} + \mu_a)(\sigma_{u\to s} + \widetilde{\phi_u} + \mu_u)} + \frac{\lambda_s \sigma_{a\to u}\sigma_{u\to s}}{(\sigma_{a\to u} + \widetilde{\phi_a} + \mu_a)(\sigma_{u\to s} + \widetilde{\phi_u} + \mu_u)(\widetilde{\phi_s} + \mu_s)} \\ &\quad + \frac{\lambda_d}{\mu_d}\left[\frac{\widetilde{\phi_a}\left(\sigma_{u\to s} + \widetilde{\phi_u} + \mu_u\right)\left(\widetilde{\phi_s} + \mu_s\right) + \widetilde{\phi_u}\sigma_{a\to u}\left(\widetilde{\phi_s} + \mu_s\right) + \widetilde{\phi_s}\sigma_{a\to u}\sigma_{u\to s}}{\left(\sigma_{a\to u} + \widetilde{\phi_a} + \mu_a\right)\left(\sigma_{u\to s} + \widetilde{\phi_u} + \mu_u\right)\left(\widetilde{\phi_s} + \mu_s\right)} - 1\right]. \end{aligned}$$

*Justification:* The typical interpretation of the $(i, j)$ −th Next Generation Matrix is the number of infections (acute individuals, members of group $i$) generated by an infected individual entering in $j$ [8]. In this manner, we may view entry $(1,4)$ of the Next Generation matrix as the approximate number of infections generated by an individual in the diagnosed stage.

From Proposition 2, we see that this entry is:

$$(NV^{-1})_{1,4} = n^T v_4 = (\lambda_a, \lambda_u, \lambda_s, \lambda_d)\begin{pmatrix} 0 \\ 0 \\ 0 \\ \dfrac{1}{\mu_d} \end{pmatrix} = \frac{\lambda_d}{\mu_d}.$$

Hence, subtracting the above from $R_t$ provides an approximation for infections generated by individuals *before* HIV diagnosis. After this point, self-testing cannot provide additional benefit towards preventing new HIV infections. Therefore, $R_{Awr}$ in a jurisdiction provides a natural indicator for the potential of HIV self-testing to reduce jurisdictional incidence.

***Theorem 1:*** At each $t$, the derivative of the susceptible population fraction:

$$\Sigma(\mathrm{t}) = \frac{susceptible\ population\ (t)}{susceptible\ population\ (t) + infected\ population\ (t)} = \frac{e}{n}$$

in (A1) obeys the bound:

$$\left|\dot{\Sigma}\right| \le \left(\frac{\Lambda}{n} + \widetilde{\lambda_a} + \mu_s - \mu_e\right)|1 - \Sigma|, \qquad (A5)$$

Where the dependence of $\Sigma$, $\dot{\Sigma}$, and $n$ on time is understood. Further, assuming $n \ge n_{min}$ for all $t$:

$$\left|\dot{\Sigma}\right| \le \left(\frac{\Lambda}{n_{min}} + \widetilde{\lambda_a} + \mu_s - \mu_e\right). \qquad (A6)$$

*Proof:* Consider (A1) and proceed by explicit calculation:

$$\begin{aligned}\frac{d}{dt}[\Sigma] &= \frac{d}{dt}\left[\frac{e}{n}\right]\\ &= \frac{\dot{e}n - \dot{n}e}{n^2}.\end{aligned}$$

$\dot{e}$ is the first equation in (A1), and $\dot{n}$ by summing the equations of (A1). Substituting into the above:

$$\begin{aligned}\left|\dot{\Sigma}\right| &= \left|\frac{\left(\Lambda - \mu_e e - \left(\widetilde{\lambda_a}a + \widetilde{\lambda_u}u + \widetilde{\lambda_s}s + \widetilde{\lambda_d}d\right)\frac{e}{n}\right)n - \left(\Lambda - (\mu_e e + \mu_a a + \mu_u u + \mu_s s + \mu_d d)\right)e}{n^2}\right|\\ &\le \left|\frac{\Lambda - \mu_e e}{n} + \frac{(\mu_e e + \mu_a a + \mu_u u + \mu_s s + \mu_d d - \Lambda)e}{n^2}\right| + \left|\frac{\left(\widetilde{\lambda_a}a + \widetilde{\lambda_u}u + \widetilde{\lambda_s}s + \widetilde{\lambda_d}d\right)e}{n^2}\right|.\end{aligned}$$

As $\left|\frac{e}{n}\right| \le 1$ by definition, $\widetilde{\lambda_a} \ge \widetilde{\lambda_{u,s,d}}$ in general, and $a + u + s + d = n - e$:

$$\begin{aligned}\left|\dot{\Sigma}\right| &\le \left|\frac{\Lambda - \mu_e e}{n} + \frac{(\mu_e e + \mu_a a + \mu_u u + \mu_s s + \mu_d d - \Lambda)e}{n^2}\right| + \left|\frac{\widetilde{\lambda_a}(n-e)}{n}\right|\\ &= \left|\frac{\Lambda - \mu_e e}{n} + \frac{(\mu_e e + \mu_a a + \mu_u u + \mu_s s + \mu_d d - \Lambda)e}{n^2}\right| + \widetilde{\lambda_a}\left|1 - \frac{e}{n}\right|.\end{aligned}$$

Again using that $\left|\frac{e}{n}\right| \le 1$, $a + u + s + d = n - e$, and additionally that $\mu_s > \mu_{e,a,u,d}$ and $\mu_e < \mu_{a,u,d,s}$ in general:

$$\begin{aligned}\left|\dot{\Sigma}\right| &\le \left|\frac{\Lambda(n-e)}{n^2} + \frac{(\mu_e e + \mu_a a + \mu_u u + \mu_s s + \mu_d d - \mu_e n)e}{n^2}\right| + \widetilde{\lambda_a}\left|1 - \frac{e}{n}\right|\\ &= \left|\frac{\Lambda(n-e)}{n^2}\right| + \left|\frac{\mu_e e + \mu_a a + \mu_u u + \mu_s s + \mu_d d - \mu_e(e + a + u + s + d)}{n}\right| + \widetilde{\lambda_a}\left|1 - \frac{e}{n}\right|\\ &= \frac{\Lambda}{n}\left|1 - \frac{e}{n}\right| + \left|\frac{(\mu_a - \mu_e)a + (\mu_u - \mu_e)u + (\mu_s - \mu_e)s + (\mu_d - \mu_e)d}{n}\right| + \widetilde{\lambda_a}\left|1 - \frac{e}{n}\right|\\ &\le \frac{\Lambda}{n}\left|1 - \frac{e}{n}\right| + \left|\frac{(\mu_s - \mu_e)(n-e)}{n}\right| + \widetilde{\lambda_a}\left|1 - \frac{e}{n}\right|\\ &= \left(\frac{\Lambda}{n} + \widetilde{\lambda_a} + \mu_s - \mu_e\right)\left|1 - \frac{e}{n}\right|,\end{aligned}$$

establishing (A5). (A6) follows immediately by noting that $0 \le \frac{e}{n} \le 1$ and $\frac{\Lambda}{n} < \frac{\Lambda}{n_{min}}$ by definition.

*Corollary of Theorem 2:* The difference between the true susceptible fraction at a time $t_n$ and the linearized susceptible fraction at time $t_0$ obeys the bound:

$$|\Sigma(t_n) - \Sigma_0| \le \left(\frac{\Lambda}{n(t^*)} + \lambda_a + \mu_s - \mu_e\right)\left(1 - \Sigma(t_n^*)\right)(t_n - t_0) \qquad (A7)$$

with $t_0 \le t_n^* < t_n$ for each $t_n$, and the bound:

$$\left(\frac{\Lambda}{n_{min}} + \lambda_a + \mu_s - \mu_e\right)(t_n - t_0) \quad (A8)$$

for all $t_n$.

*Proof:* Since $e$ and $n$ are both at least once differentiable and nonzero by definition, $\Sigma = \frac{e}{n}$ may be expanded in a Taylor series around $\Sigma_0$ for each $t_n$ such that:

$$\Sigma(t_n) = \Sigma_0 + \dot{\Sigma}(t_n^*)(t_n - t_0),$$

for some $t_n^*$ in $[t_0, t_n]$. From (A5):

$$\begin{aligned} |\Sigma(t_n) - \Sigma_0| &= \left|\dot{\Sigma}(t_n^*)\right||(t_n - t_0)| \\ &\le \left(\frac{\Lambda}{n(t_n^*)} + \widetilde{\lambda_a} + \mu_s - \mu_e\right)|1 - \Sigma(t_n^*)|(t_n - t_0), \end{aligned}$$

Establishing (A7). (A8) follows by an identical argument using (A6) rather than (A5).

Note that the bound (A8) is not a particularly sharp bound and can likely be improved. Nevertheless, it provides a uniform bound.

In practice, however, (A7) provides an effective estimate. Since $\Sigma_0 \approx 1$ in general, and, by continuity of $\Sigma$, we then expect $|1 - \Sigma(t)| \approx 0$ for $t$ close to $t_0$, ensuring (A7) remains sharp for $t$ not excessively large.

Note that $\widetilde{\lambda_a}$, while convenient insofar as providing an upper bound for analysis, is generally much higher than the overall transmission rate, and hence not representative of its true value in practice. Additionally, $\frac{\Lambda}{n}$ is quite small in general, as $n$ is several orders of magnitude larger than $\Lambda$.

Note further that (A7) is thus easy to interpret: it states that our linearized approximation of the system (A1) remains valid provided:

- New entries to the eligible (susceptible) population are small compared to the overall population.
- Mortality and transmission remain low.
- The prevalence of HIV in the population is low.

The scenarios examined in the current document satisfy these conditions, and so we expect the linearized system (1) in the main text to provide a good approximation to (A1) throughout our time horizon.

**Supplement B**

*Compartment model parameterization*

This section outlines how we defined the disease-stage specific transmission, mortality, and testing parameters mentioned in the methods section. We used measured data from ATLAS (2017-2019) for simulations in the main text and PATH outputs from the baseline PATH simulation for validation in supplement C. In both cases, we considered the full-population level transmission rate defined as:

$$\bar{\lambda} = \frac{Incidence}{Prevalence} \qquad \text{(B1).}$$

averaged over the time-period. In the main text, this information, along with other relevant surveillance data, came from ATLAS database [16]. We assume for any transmission level, the $\lambda_i$ are related by multiplicative factors $\alpha_i$. These factors were derived from [21] and reported in Table 1 in the main text.

The individual $\lambda_i$ are given by:

$$\lambda_u = \frac{\bar{\lambda}}{\alpha_a P(A) + P(U) + \alpha_s P(S) + \alpha_{noCare} P(noCare) + \alpha_{ART} P(ART) + \alpha_{VLS} P(VLS)},$$
$$\lambda_a = \alpha_a \lambda_u,$$
$$\lambda_s = \alpha_s \lambda_u,$$
$$\lambda_{noCare} = \alpha_{noCare} \lambda_u,$$
$$\lambda_{ART} = \alpha_{ART} \lambda_u$$
$$\lambda_{VLS} = \alpha_{VLS} \lambda_{VLS}.$$

$P(ART), P(noCare), P(VLS)$ were available from our measured data [16], [21].

$P(A),\ P(U), P(S)$, referring to the probabilities of being in the acute, chronic unaware, and AIDS unaware states, respectively, were obtained by first considering the probability that an individual is unaware of their HIV status $P(\bar{D})$, available from data:

$$P(\bar{D}) = 1 - P(D).$$

$P(A|\bar{D})$ and $P(S|\bar{D})$ were determined through calibration. Accordingly:

$$P(A) = P(A|\bar{D}) \times P(\bar{D}), \qquad P(S) = P(S|\bar{D}) \times P(\bar{D})$$

and

$$P(U) = P(\bar{D}) - P(A) - P(S).$$

Mortality rates were obtained similarly, with

$$\bar{\mu} = \frac{PWHDeaths}{Prevalence}, \qquad \text{(B2)}$$

and the data coming from ATLAS. This yields:

$$\mu_a = \frac{\bar{\mu}}{P(A) + \beta_u P(U) + \beta_s P(S) + \beta_{noCare} P(noCare) + \beta_{ART} P(ART) + \beta_{VLS} P(VLS)},$$

$\mu_u = \beta_u \mu_a,$
$\mu_s = \beta_s \mu_a,$
$\mu_{noCare} = \beta_{noCare} \mu_a,$
$\mu_{ART} = \beta_{ART} \mu_a$
$\mu_{VLS} = \beta_{VLS} \mu_{VLS}.$

The $\beta_i$ also came from [21] and are reported in Table 1 in the main text.

Testing rates (shown in Table 1 in the main text) were approximated as follows.

First, using surveillance data, we define the following quantity:

$$\overline{\phi} \approx \frac{annual\ new\ diagnoses}{PWH\ unaware\ in\ previous\ year + annual\ new\ infection}. \qquad \text{(B3)}$$

This gives us the rate at which undiagnosed PWH receive a diagnosis.

Next, assume that undiagnosed PWH with chronic infection test at some unknown rate $\phi_u$, and further, that undiagnosed PWH with acute infection and AIDS test at $\nu_a$ and $\nu_s$ this rate, respectively. Then $\phi_a = \nu_a \phi_u, \phi_s = \nu_s \phi_u$.

Let $\overline{D}$ refer to the probability of being undiagnosed. Then the probability of having undiagnosed acute and undiagnosed late-stage infection (AIDS) are denoted:

$$P(A|\overline{D}), \qquad P(S|\overline{D}).$$

The probability of an undiagnosed PWH having a chronic-stage infection is thus obtained as:

$$1 - P(A|\overline{D}) - P(S|\overline{D}).$$

Finally, we account for test sensitivity, assumed to vary by stage: $\kappa_{a,s,u}$.

The overall testing rate among undiagnosed PWH is:

$$\phi_{undiag} = P(A|\overline{D})\phi_a + P(S|\overline{D})\phi_s + \left(1 - P(A|\overline{D}) - P(S|\overline{D})\right)\phi_u,$$

which reduces to:

$$\phi_{undiag} = \left(P(A|\overline{D})\nu_a + P(S|\overline{D})\nu_s + \left(1 - P(A|\overline{D}) - P(S|\overline{D})\right)\right)\phi_u.$$

Denote stage-specific test sensitivities as $\kappa_i$. Assuming every positive test among an undiagnosed PWH leads to a diagnosis, then we have, approximately:

$$\overline{\phi} = \left(\kappa_a P(A|\overline{D})\nu_a + \kappa_s P(S|\overline{D})\nu_s + \kappa_u\left(1 - P(A|\overline{D}) - P(S|\overline{D})\right)\right)\phi_u,$$

implying:

$$\phi_u = \frac{\overline{\phi}}{\kappa_a \nu_a \alpha_a P(A|\overline{D}) + \nu_s \kappa_s \alpha_s P(S|\overline{D}) + \alpha_u \kappa_u\left(1 - P(A|\overline{D}) - P(S|\overline{D})\right)}.$$

Allowing us to recover as $\phi_{a,s}$ as:

$$\phi_a = \nu_a \phi_u, \quad \phi_s = \nu_s \phi_u.$$

For each of the 38 counties studied, we evaluate the model parameterization procedure by comparing the surveillance data against simulation data over the parameterization years (2017-19) for the following:

1. Incidence rate per month $\overline{\lambda}$ (Eq. B1);
2. PWH mortality rate per year $\overline{\mu}$ (Eq. B2);
3. % of PWH aware of status;
4. PWH testing rate per month $\overline{\phi}$ (Eq. B3).

The results are provided in Table B1 and Figures B1-B5.

**Table B1:** Agreement between surveillance and simulation for key HIV indicators for each jurisdiction. Definitions for $\overline{\lambda}, \overline{\mu}, \overline{\phi}$ are given in equations B1-B3, respectively. For surveillance quantities, the three-year average is reported, with the three-year range provided in parentheses.

| Jurisdiction | Surv. $\overline{\lambda}$, 2017-19 | Sim. $\overline{\lambda}$, 2017-19 | Surv. $\overline{\mu}$, 2017-19 | Sim. $\overline{\mu}$, 2017-19 | Surv. % aware, 2017-19 | Sim. % aware, 2017-19 | Surv. $\overline{\phi}$, 2017-19 | Sim. $\overline{\phi}$, 2017-19 |
|---|---|---|---|---|---|---|---|---|
| Alameda County, CA | 0.028 (0.028-0.029) | 0.028 | 0.009 (0.009-0.010) | 0.009 | 87.2 (86.8-87.6) | 88.0 | 0.016 (0.016-0.017) | 0.016 |
| Los Angeles County, CA | 0.026 (0.025-0.028) | 0.027 | 0.010 (0.009-0.010) | 0.010 | 88.6 (88.1-89.0) | 89.2 | 0.018 (0.017-0.018) | 0.018 |
| Riverside County, CA | 0.027 (0.025-0.029) | 0.028 | 0.015 (0.013-0.017) | 0.016 | 88.8 (88.4-89.3) | 88.3 | 0.016 (0.015-0.017) | 0.016 |
| Sacramento County, CA | 0.042 (0.040-0.044) | 0.048 | 0.012 (0.010-0.015) | 0.013 | 83.8 (83.5-84.1) | 80.2 | 0.014 (0.013-0.015) | 0.014 |
| San Bernardino County, CA | 0.050 (0.046-0.053) | 0.051 | 0.013 (0.010-0.016) | 0.013 | 78.4 (77.4-79.4) | 79.4 | 0.015 (0.015-0.016) | 0.015 |
| San Diego County, CA | 0.028 (0.027-0.029) | 0.029 | 0.009 (0.008-0.009) | 0.009 | 86.3 (86.3-86.4) | 86.2 | 0.013 (0.013-0.014) | 0.014 |
| San Francisco County, CA | 0.012 (0.010-0.013) | 0.012 | 0.008 (0.001-0.015) | 0.007 | 96.2 (95.6-96.9) | 96.6 | 0.031 (0.028-0.035) | 0.028 |
| Broward County, FL | 0.027 (0.024-0.030) | 0.027 | 0.013 (0.012-0.014) | 0.014 | 89.1 (88.8-89.4) | 89.0 | 0.017 (0.016-0.018) | 0.017 |
| Duval County, FL | 0.034 (0.030-0.038) | 0.033 | 0.018 (0.017-0.018) | 0.017 | 83.2 (82.4-84.1) | 86.2 | 0.016 (0.015-0.017) | 0.016 |
| Hillsborough County, FL | 0.034 (0.029-0.038) | 0.033 | 0.015 (0.013-0.016) | 0.015 | 84.3 (83.8-84.8) | 85.6 | 0.015 (0.014-0.016) | 0.015 |
| Miami-Dade County, FL | 0.035 (0.033-0.037) | 0.038 | 0.013 (0.012-0.013) | 0.013 | 87.1 (86.8-87.4) | 86.1 | 0.019 (0.019-0.019) | 0.018 |
| Palm Beach County, FL | 0.028 (0.026-0.029) | 0.028 | 0.018 (0.017-0.019) | 0.018 | 86.7 (86.6-86.8) | 87.5 | 0.015 (0.013-0.016) | 0.015 |
| Pinellas County, FL | 0.029 (0.027-0.031) | 0.027 | 0.020 (0.018-0.021) | 0.019 | 86.7 (85.8-87.6) | 88.9 | 0.017 (0.015-0.018) | 0.017 |
| Cobb County, GA | 0.046 (0.041-0.051) | 0.050 | 0.007 (0.007-0.008) | 0.008 | 81.5 (80.9-82.1) | 79.8 | 0.015 (0.014-0.016) | 0.015 |
| Dekalb County, GA | 0.035 (0.031-0.039) | 0.037 | 0.009 (0.009-0.009) | 0.009 | 84.6 (84.1-85.1) | 84.5 | 0.016 (0.015-0.016) | 0.016 |
| Fulton County, GA | 0.032 (0.029-0.035) | 0.034 | 0.009 (0.009-0.010) | 0.010 | 85.4 (84.9-85.9) | 85.4 | 0.015 (0.014-0.016) | 0.016 |
| Gwinnett County, GA | 0.047 (0.043-0.051) | 0.048 | 0.006 (0.005-0.007) | 0.007 | 80.4 (79.4-81.5) | 81.7 | 0.017 (0.015-0.020) | 0.017 |
| Marion County, IN | 0.041 (0.041-0.042) | 0.045 | 0.010 (0.009-0.012) | 0.010 | 84.5 (84.2-84.7) | 82.7 | 0.017 (0.015-0.018) | 0.016 |
| East Baton Rouge Parish, LA | 0.037 (0.035-0.040) | 0.039 | 0.019 (0.018-0.020) | 0.019 | 85.0 (84.9-85.1) | 84.6 | 0.016 (0.015-0.017) | 0.016 |
| Orleans Parish, LA | 0.026 (0.024-0.028) | 0.027 | 0.014 (0.013-0.016) | 0.014 | 89.1 (88.5-89.8) | 90.0 | 0.020 (0.019-0.020) | 0.019 |
| Montgomery County, MD | 0.022 (0.021-0.024) | 0.022 | 0.006 (0.004-0.007) | 0.005 | 89.6 (88.8-90.4) | 91.8 | 0.022 (0.020-0.024) | 0.020 |
| Prince George's County, MD | 0.026 (0.023-0.029) | 0.025 | 0.006 (0.005-0.007) | 0.006 | 86.8 (86.0-87.5) | 89.5 | 0.017 (0.017-0.018) | 0.017 |
| Baltimore City, MD | 0.016 (0.012-0.019) | 0.016 | 0.017 (0.010-0.024) | 0.018 | 92.1 (91.8-92.3) | 93.3 | 0.017 (0.017-0.018) | 0.017 |
| Wayne County, MI | 0.036 (0.035-0.037) | 0.036 | 0.015 (0.011-0.019) | 0.016 | 85.6 (84.8-86.3) | 86.4 | 0.018 (0.018-0.019) | 0.018 |
| Bronx County, NY | 0.014 (0.013-0.014) | 0.014 | 0.016 (0.015-0.017) | 0.016 | 94.0 (93.6-94.4) | 94.7 | 0.020 (0.018-0.021) | 0.019 |
| Kings County, NY | 0.016 (0.014-0.017) | 0.016 | 0.014 (0.014-0.015) | 0.014 | 92.7 (92.4-93.0) | 93.6 | 0.019 (0.017-0.021) | 0.018 |
| New York County, NY | 0.011 (0.010-0.012) | 0.011 | 0.010 (0.009-0.011) | 0.010 | 94.4 (94.2-94.6) | 95.0 | 0.016 (0.016-0.017) | 0.016 |
| Queens County, NY | 0.017 (0.015-0.020) | 0.017 | 0.008 (0.007-0.010) | 0.008 | 91.0 (90.4-91.7) | 92.8 | 0.018 (0.018-0.019) | 0.018 |
| Cuyahoga County, OH | 0.024 (0.024-0.025) | 0.024 | 0.013 (0.010-0.016) | 0.012 | 86.6 (86.0-87.3) | 89.1 | 0.015 (0.014-0.016) | 0.015 |
| Franklin County, OH | 0.035 (0.035-0.036) | 0.036 | 0.011 (0.008-0.013) | 0.011 | 85.4 (85.0-85.8) | 85.6 | 0.017 (0.016-0.018) | 0.017 |
| Hamilton County, OH | 0.039 (0.036-0.042) | 0.038 | 0.014 (0.013-0.015) | 0.014 | 81.5 (80.1-83.0) | 85.6 | 0.019 (0.018-0.019) | 0.018 |
| Philadelphia County, PA | 0.022 (0.021-0.023) | 0.023 | 0.014 (0.013-0.015) | 0.014 | 91.9 (91.7-92.1) | 91.7 | 0.021 (0.019-0.022) | 0.019 |
| Shelby County, TN | 0.036 (0.034-0.038) | 0.036 | 0.015 (0.014-0.016) | 0.015 | 85.7 (85.2-86.2) | 85.9 | 0.018 (0.016-0.019) | 0.017 |
| Dallas County, TX | 0.035 (0.030-0.040) | 0.037 | 0.007 (0.003-0.010) | 0.008 | 83.8 (83.2-84.4) | 84.4 | 0.016 (0.015-0.016) | 0.016 |
| Harris County, TX | 0.037 (0.036-0.038) | 0.038 | 0.013 (0.012-0.014) | 0.013 | 83.7 (83.0-84.4) | 84.5 | 0.016 (0.015-0.017) | 0.016 |
| Tarrant County, TX | 0.047 (0.041-0.053) | 0.049 | 0.012 (0.012-0.013) | 0.013 | 82.0 (81.7-82.3) | 80.8 | 0.016 (0.016-0.017) | 0.016 |
| Travis County, TX | 0.035 (0.030-0.041) | 0.037 | 0.010 (0.009-0.011) | 0.010 | 83.5 (83.3-83.7) | 83.2 | 0.014 (0.013-0.016) | 0.014 |
| King County, WA | 0.030 (0.027-0.033) | 0.029 | 0.011 (0.010-0.011) | 0.010 | 87.5 (87.1-88.0) | 88.3 | 0.017 (0.015-0.019) | 0.017 |

**Figure B1:** Surveillance vs. simulation, jurisdictions 1-8

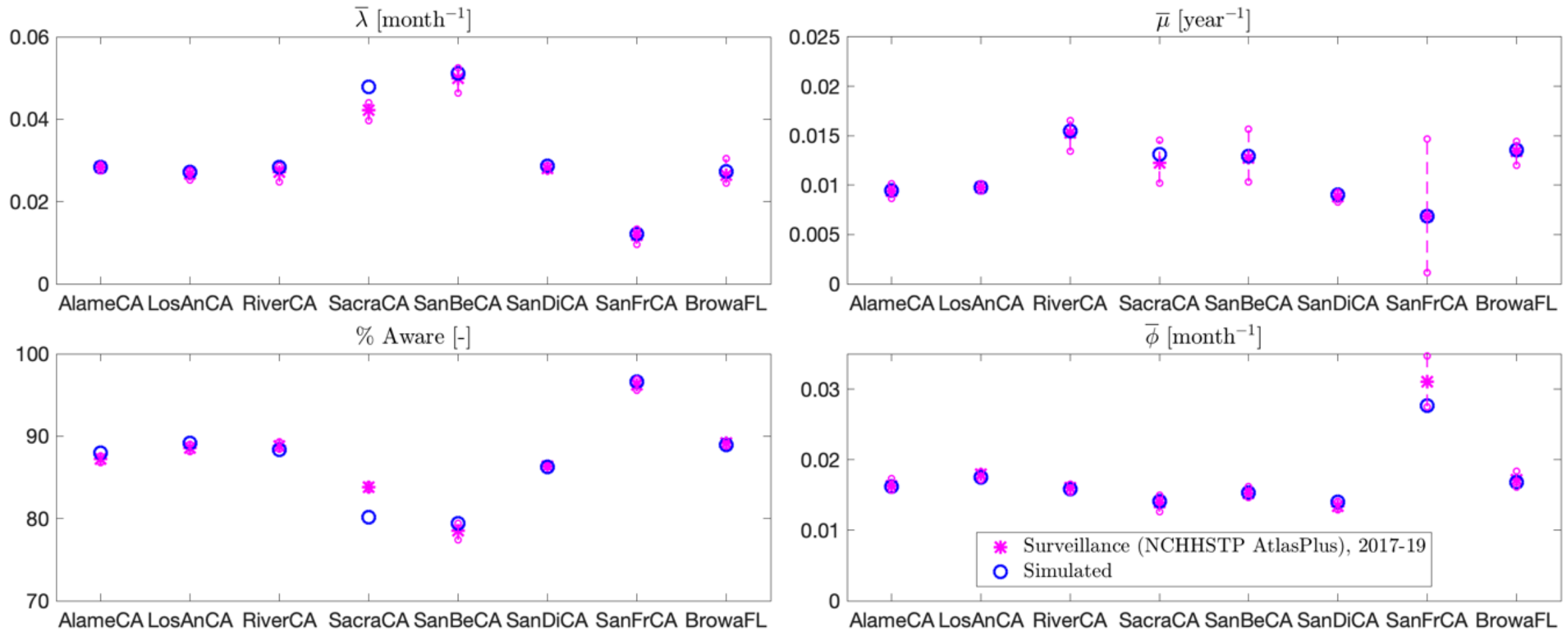

**Figure B2:** Surveillance vs. simulation, jurisdictions 9-17

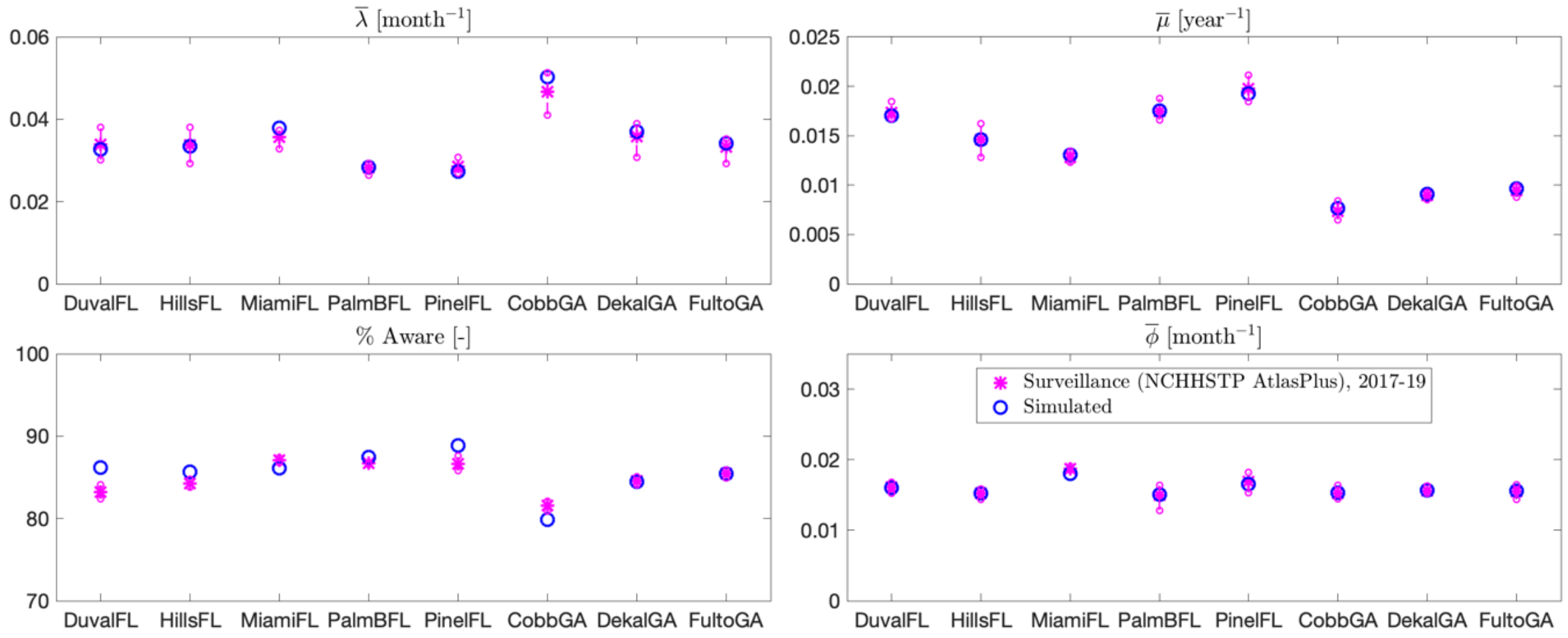

**Figure B3:** Surveillance vs. simulation, jurisdictions 18-25

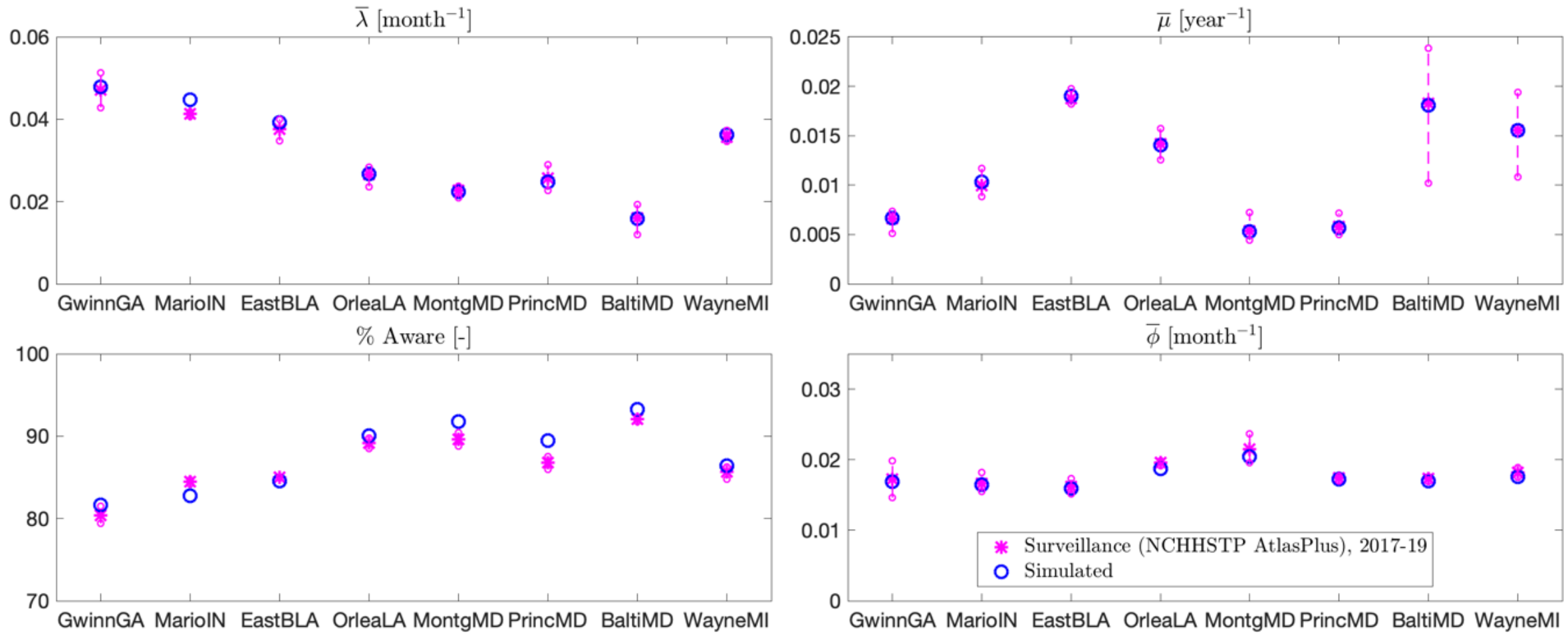

**Figure B4:** Surveillance vs. simulation, jurisdictions 26-33

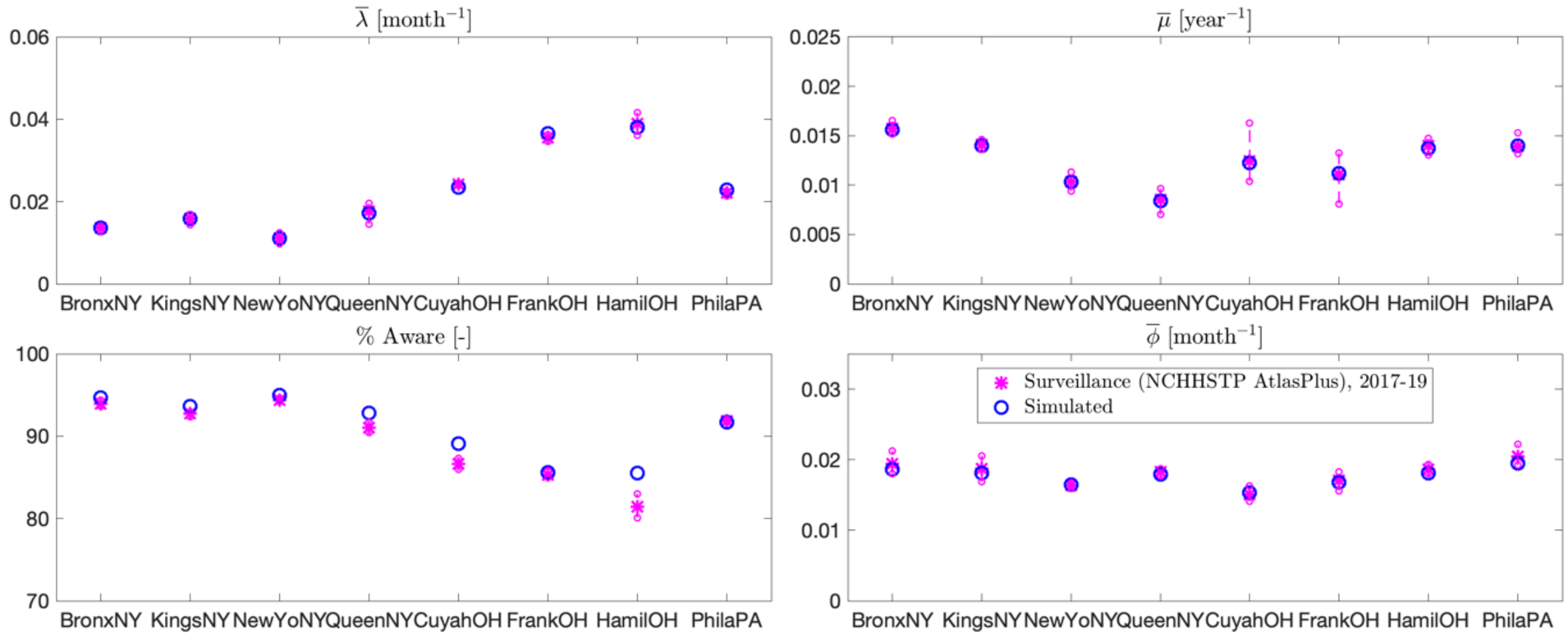

**Figure B5:** Surveillance vs. simulation, jurisdictions 34-38

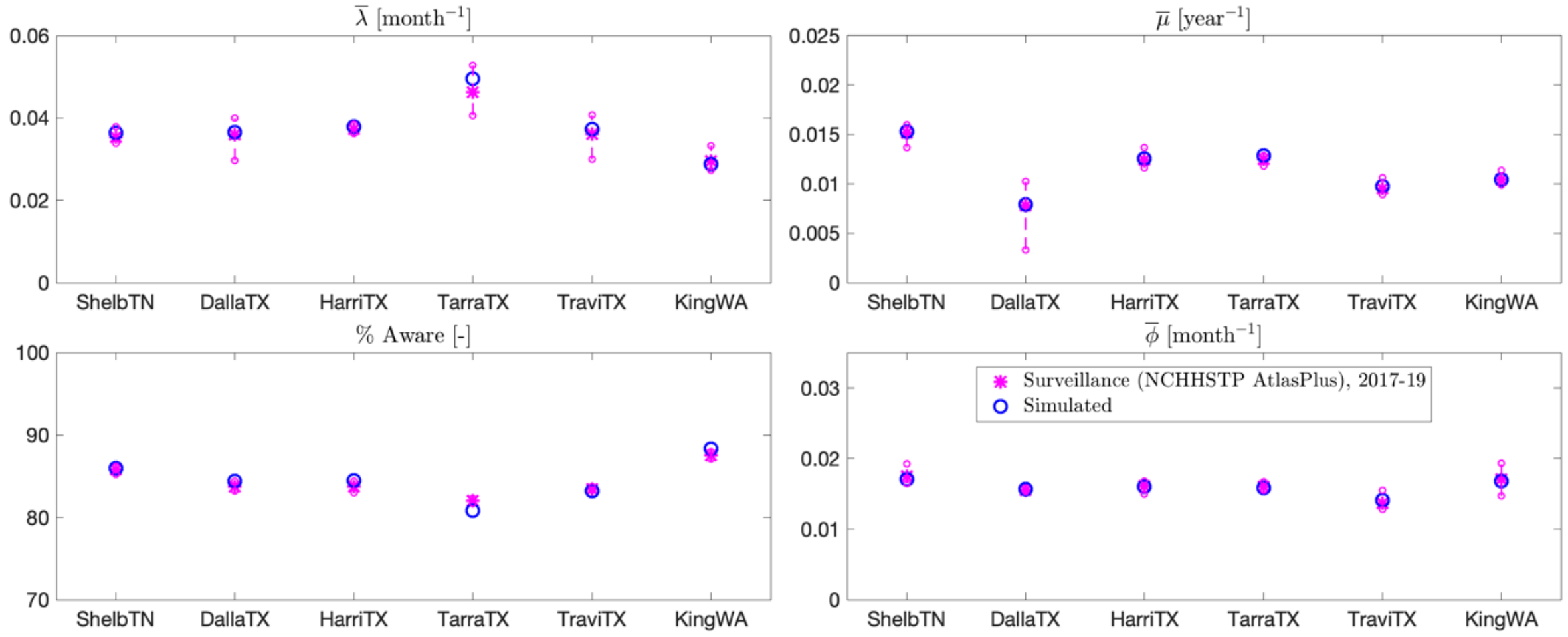